\documentclass[11pt, letterpaper]{article}
\usepackage[brazil,english]{babel}
\usepackage[latin1]{inputenc}
\usepackage{fullpage}
\usepackage[dvips]{graphicx}
\usepackage{amssymb, amsmath}
\usepackage{color}
\usepackage{footnote}
\newtheorem{theorem}{Theorem}[section]

\newtheorem{definition}[theorem]{Definition} 
\newtheorem{example}[theorem]{Example}  
\newtheorem{remark}[theorem]{Remark}  

\def\beginproof{\noindent{\em Proof}.~}
\def\endproof{{\ \hfill\hbox{\fbox{}}\parfillskip 0pt}\par}
\def\beginsol{\noindent{\em Resolution}.~}
\def\endsol{{\ \hfill\hbox{\fbox{}}\parfillskip 0pt}\par}

\renewcommand{\thefootnote}{\arabic{ftnote}}
\newcommand{\R}{\mathbb{R}}
\title{Calculus, constrained minimization and Lagrange multipliers: \\ Is the 
optimal critical point a local minimizer?
\footnote{This short note is for those who study or teach calculus.}}
\author{Ademir Alves Ribeiro\footnotemark[4] \and 
Jos\'e Renato Ramos Barbosa\footnotemark[4]}
\begin{document}
\maketitle
\renewcommand{\thefootnote}{\fnsymbol{footnote}}
\footnotetext[4]{Department of Mathematics, Federal University of Paran\'a, 
Brazil (email: {ademir.ribeiro@ufpr.br, jrrb@ufpr.br}). The first author is 
supported by CNPq, Brazil, Grant 309437/2016-4.}

\renewcommand{\thefootnote}{\arabic{ftnote}}
\begin{abstract}
In this short note, we discuss how the optimality conditions for the problem of
minimizing a multivariate function subject to equality constraints have been
dealt with in undergraduate Calculus. We are particularly interested in
the $2$ or $3$-dimensional cases, which are the most common cases in Calculus 
courses. Besides giving sufficient conditions to a critical point to be a local
minimizer, we also present and discuss counterexamples to some statements 
encountered in the undergraduate literature on Lagrange Multipliers, such as 
``among the critical points, the ones which have the smallest image (under the 
function) are minimizers" or ``a single critical point (which is a local minimizer) is a global 
minimizer".
\end{abstract}


\thispagestyle{plain}

\section{Introduction}
In spite of being a strategy for finding the local maxima and/or minima of a
function subject to constraints, the Lagrange Multiplier Method (LMM), 
particularly when it is used for solving undergraduate Optimization problems, is
tipically used as a systematic procedure for identifying global extrema. In 
doing so, some undergraduate textbooks on the subject show the statements 
related to the LMM (and/or the worked problems based on it) without complete 
hypotheses or carelessly written in an imprecise manner and with 
oversimplification (as an attempt to make the method more palatable). 
Specifically, the validity of the LMM, when one is looking for global extrema, 
depends on the existence of those extrema. This basic assumption has to be 
satisfied beforehand. Otherwise, even if one succeeded in obtaining 
local extrema from the critical points determined by the method, it might not be
possible to get the global extrema from the local ones. However, it is not 
unusual to find books or academic homepages where, right after the local extrema
are found, they are promptly evaluated and the ones with the smallest 
(respectively, the greatest) image under the function are elected the global 
minima (respectively, maxima). This is done even without having been previously 
established the existence of the global extrema. The problem gets worse when 
there is just one critical point and one has to determine the maximum/minimum 
value from it without having a local criterion to be used in conjuction with the
LMM. We present such a criterion here and propose the adoption of it in Calculus
courses. At least, such a procedure makes the student completely sure that the 
optimal local values can be determined from the critical points. Concerning the
global aspect, well, this is a whole different story. Many times, without 
compactness, the proof that there are global extrema is something out of the 
scope of Calculus courses and depends on the specific problem that we are 
dealing with. On the other hand, even in the very few books which correctly 
state the LMM, such as the excellent \cite{marsden}, emphasizing only the local
aspect of the method, when working on problems with lack of compactness, it is 
assumed, for instance, the existence of ``a box of largest possible volume" 
among all rectangular boxes with a fixed surface area. Then, right after finding
a unique critical point via the LMM, the solution ends with a conclusion like 
this one: ``This (cubical) shape must therefore maximize the volume, assuming 
there is a box of maximum volume." In this paper we work on a similar problem, 
showing the existence of the global extrema via a nontrivial reasoning for 
Calculus courses. Furthermore, we emphasize that a criterion to determine if a 
critical point (obtained by the LMM) is a local maximum/minimum would help 
enormously in problems like those ones.

\section{Preliminaries}
\label{sec_intro}
We consider here the equality constrained optimization problem of the form
\begin{equation}
\label{eq_prob}
\begin{array}{cl}
\displaystyle{\rm minimize } & f(x) \\
{\rm subject\ to} & g(x)=0,
\end{array}
\end{equation}
where $f:D\to\R$ and $g:D\to\R^m$ are twice continuously differentiable 
functions defined on the open set $D\subset\R^n$. The set
\begin{equation}
\label{feas_set}
\Omega=\{x\in D\mid g(x)=0\}
\end{equation}
is called {\em feasible set} of the problem (\ref{eq_prob}). Recall the well 
known definition of a minimizer:
\begin{definition}
\label{def:min}
A point $x^*\in\Omega$ is said to be a local minimizer of the problem 
(\ref{eq_prob}) if there exists $\delta>0$ such that $f(x^*)\leq f(x)$ for all 
$x\in\Omega\cap B(x^*,\delta)$. If $f(x^*)\leq f(x)$ for all $x\in\Omega$, the 
point $x^*$ is called a global minimizer. 
\end{definition}

\begin{remark}
\label{rmk:max}
{\em There is no loss of generality in considering only minimization problems 
since if we want to maximize a function $f$, we can equivalently minimize $-f$. 
So, the definitions and results can be easily rewritten.}
\end{remark}

In the following example we can directly verify that a point is a minimizer.
Nevertheless, this is not always the case and we normally need other tools to 
find minimizers. It should be pointed out that, in spite of focusing on 
constrained problems, some examples can be better understood and/or visualized 
if we disregard the constraint. In fact, we can transform an unconstrained 
problem into an equivalent constrained one by introducing an artificial 
variable: 
$$
\begin{array}{cl}
\displaystyle{\rm minimize } & f(x)
\end{array}
\quad\mbox{is equivalent to}\quad
\begin{array}{cl}
\displaystyle{\rm minimize } & \varphi(x,y):=f(x) \\
{\rm subject\ to} & g(x,y):=y=0,
\end{array}
$$

\begin{figure}[htbp]
\centering
\includegraphics[scale=0.5]{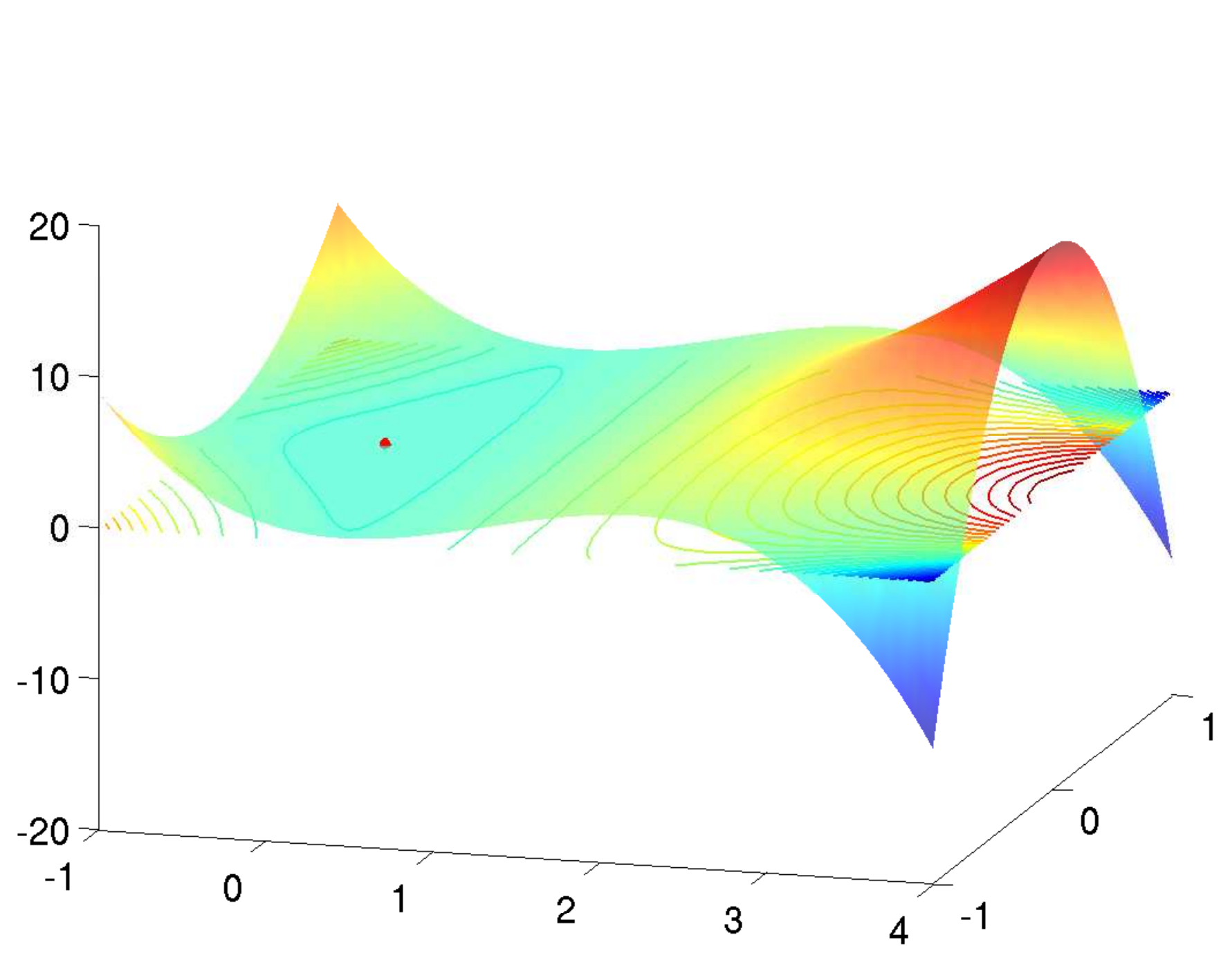}
\caption{Graph of the function $f$, given in Example \ref{ex:local}, showing that 
the local minimizer is not a global minimizer.}
\label{fig_ctex1}
\end{figure}

\begin{example}
\label{ex:local}
Consider the function $f:\R^2\to\R$ given by $f(x)=x_1^2+x_2^2(1-x_1)^3$.
The point $x^*=0$ is a local minimizer since $f(x^*)=0\leq f(x)$ for all 
$x\in B(x^*,1)$. Moreover, this point is not a global minimizer because 
$f(4,1)=-11$. See Figure \ref{fig_ctex1}.
\end{example}

A well known condition that ensures the existence of a global minimizer is the 
compactness of the set. 

\begin{theorem}
\label{th:compact}
Let $L\subset\Omega$ a compact set. Then $f$ has a global minimizer in $L$, that
is, there exists $x^*\in L$ such that $f(x^*)\leq f(x)$ for all $x\in L$.
\end{theorem}

\noindent Unfortunately, there are many situations where the above result cannot
be directly applied since the underlying set might not be compact. Even so, we 
may still guarantee the existence of a global minimizer, as discussed in 
the upcoming Example \ref{ex:min_area}. 

\section{Necessary optimality conditions: Lagrange multipliers}
\label{sec_lagrange}
In this section we present the necessary conditions that must be satisfied by every local 
minimizer. We also point out two misunderstandings that sometimes arise in calculus courses.
\begin{definition}
\label{def:stat}
A point $x^*\in\R^n$ is said to be critical (or stationary) for the problem (\ref{eq_prob}) 
if there exists a vector $\lambda^*\in\R^m$  such that 
\begin{subequations}
\begin{align}
\nabla f(x^*)+\sum_{i=1}^{m}\lambda_i^*\nabla g_i(x^*)=0, & \label{kkt_grad} \\[-3pt]
g(x^*)=0 & \label{kkt_feas}.
\end{align}
\end{subequations}
\end{definition}
The components of $\lambda^*$ are the {\em Lagrange multipliers} associated with
the constraints.

The next result is a classical one and it is used to find possible candidates 
for the optimal solutions. (See, for instance, \cite{fitzpatrick} for a version
of such theorem.)
\begin{theorem}
\label{th:kkt}
Suppose that $x^*\in\R^n$ is a local minimizer for the problem (\ref{eq_prob}) 
and the gradients $\nabla g_i(x^*)$ are linearly independent, $i=1,\ldots,m$. 
Then $x^*$ is a critical point for this problem. 
\end{theorem}

As previously mentioned, a very common exercise in undergraduate Calculus is the
problem of minimizing the area of a box, without lid, subject to a constant 
volume. The issue here rely on the fact that, normally, the solution is
not accompanied by a mathematical argument explaining that the box with minimum
area does exist and/or with a justification as to why the critical point 
obtained (via the equations (\ref{kkt_grad})--(\ref{kkt_feas})) is the global 
minimizer of the problem. Let us discuss these issues more precisely in the next
example.

\begin{example}
\label{ex:min_area}
Consider $D=\{x\in\R^3\mid x_1>0, x_2>0, x_3>0\}$ and the functions $f,g:D\to\R$ defined 
by $f(x)=x_1x_2+2x_1x_3+2x_2x_3$ and $g(x)=x_1x_2x_3-1$. Show that the problem (\ref{eq_prob}) 
has a (unique) global solution and find it using the Lagrange method.
\end{example}
\beginsol
As defined before, let $\Omega=\{x\in D\mid g(x)=0\}$ be the feasible set of the problem. 
We claim that if $x\in\Omega$ and $f(x)\leq 5$, then $\frac{1}{5}\leq x_i\leq 5$, for 
$i=1,2,3$. Indeed, if $x_1<\frac{1}{5}$ or $x_2<\frac{1}{5}$, then 
$$
f(x)=x_1x_2+2x_1x_3+\dfrac{2}{x_1}>\dfrac{1}{x_1}>5\quad\mbox{or}\quad
f(x)=x_1x_2+\dfrac{2}{x_2}+2x_2x_3>\dfrac{1}{x_2}>5, 
$$
respectively. If $x_3<\frac{1}{5}$, then 
$f(x)=\dfrac{1}{x_3}+2x_1x_3+2x_2x_3>\dfrac{1}{x_3}>5$. Furthermore, if $x_1>5$ or $x_2>5$, 
$$
f(x)=x_1x_2+\dfrac{2(x_1+x_2)}{x_1x_2}>x_1x_2+\dfrac{10}{x_1x_2}=
\dfrac{1}{x_1x_2}\left(\left(x_1x_2-\dfrac{5}{2}\right)^2+\dfrac{15}{4}\right)+5>5.
$$
If $x_3>5$, then 
$$
f(x)=\dfrac{x_2+2x_3}{x_2x_3}+2x_2x_3>\dfrac{10}{x_2x_3}+2x_2x_3=
\dfrac{2}{x_2x_3}\left(\left(x_1x_2-\dfrac{5}{4}\right)^2+\dfrac{55}{16}\right)+5>5,
$$
proving the claim. Now, consider the set $L=\{x\in\Omega\mid f(x)\leq 5\}$. If 
$(x^k)\subset L$ is such that $x^k\to\bar x$, then $\bar x\in D$, $g(\bar x)=0$ and 
$f(\bar x)\leq 5$, which means that $L$ is closed. It is also bounded in view of
the claim. Thus, Theorem \ref{th:compact} ensures that there exists $x^*\in L$ 
such that $f(x^*)\leq f(x)$ for all $x\in L$. This point is in fact a global 
minimizer in $\Omega$, because if $x\in\Omega\setminus L$, we have 
$f(x)>5\geq f(x^*)$.\\
\noindent Now, applying Theorem \ref{th:kkt}, we conclude that $x^*$ must be 
solution of the equations 
$$
\left(\begin{array}{c} x_2+2x_3 \\ x_1+2x_3 \\ 2x_1+2x_2 \end{array}\right)+\lambda
\left(\begin{array}{c} x_2x_3 \\ x_1x_3 \\ x_1x_2 \end{array}\right)=
\left(\begin{array}{c} 0 \\ 0 \\ 0 \end{array}\right)\quad\mbox{and}\quad x_1x_2x_3=1.
$$
Since this system has a unique solution, namely, 
$\dfrac{1}{2}\left(\begin{array}{c} 2\sqrt[3]{2} \\ 2\sqrt[3]{2} \\ \sqrt[3]{2}
\end{array}\right)$ with $\lambda=-2\sqrt[3]{4}$, it follows that $x^*$ is 
exactly this point.
\endsol
\noindent Concerning the previous example, it is worth mentioning that,
on the one hand, its reasoning is out of the scope of undergraduate Calculus 
textbooks. On the other hand, it is also a remainder that, sometimes, it is not 
trivial to establish the existence of a global minimum.

The next example is a reformulation of the unconstrained problem given in 
Example \ref{ex:local} as an equivalent constrained problem, obtained by 
introducing an artificial variable. 

\begin{example}
\label{ex:local_only1}
Consider the functions $f,g:\R^3\to\R$ given by $f(x)=x_1^2+x_2^2(1-x_1)^3$ and $g(x)=x_3$. 
It is easy to verify that the point $x^*=0$ and the multiplier $\lambda^*=0$ satisfy the 
conditions (\ref{kkt_grad})--(\ref{kkt_feas}). In fact, it can be proved that this point 
is the only critical point of this example. Moreover, $x^*$ is a local minimizer since 
$f(x^*)=0\leq f(x)$ for all $x\in B(x^*,1)$. However, this point is not a global minimizer 
because $f(4,1,0)=-11$. 
\end{example}

\begin{remark}
\label{rmk:ctex1}
{\em Note that the previous example also answers negatively the question 
{\bf ``If a function has a single critical point which is a local minimizer, is
this point a global minimizer?"}, which sometimes takes place (in some Calculus
courses and academic homepages) and is responded incorrectly with a {\bf``yes"}.
This probably occurs since, for functions of one variable, the result holds 
as the next theorem states.}
\end{remark}

\begin{figure}[htbp]
\centering
\includegraphics[scale=0.385]{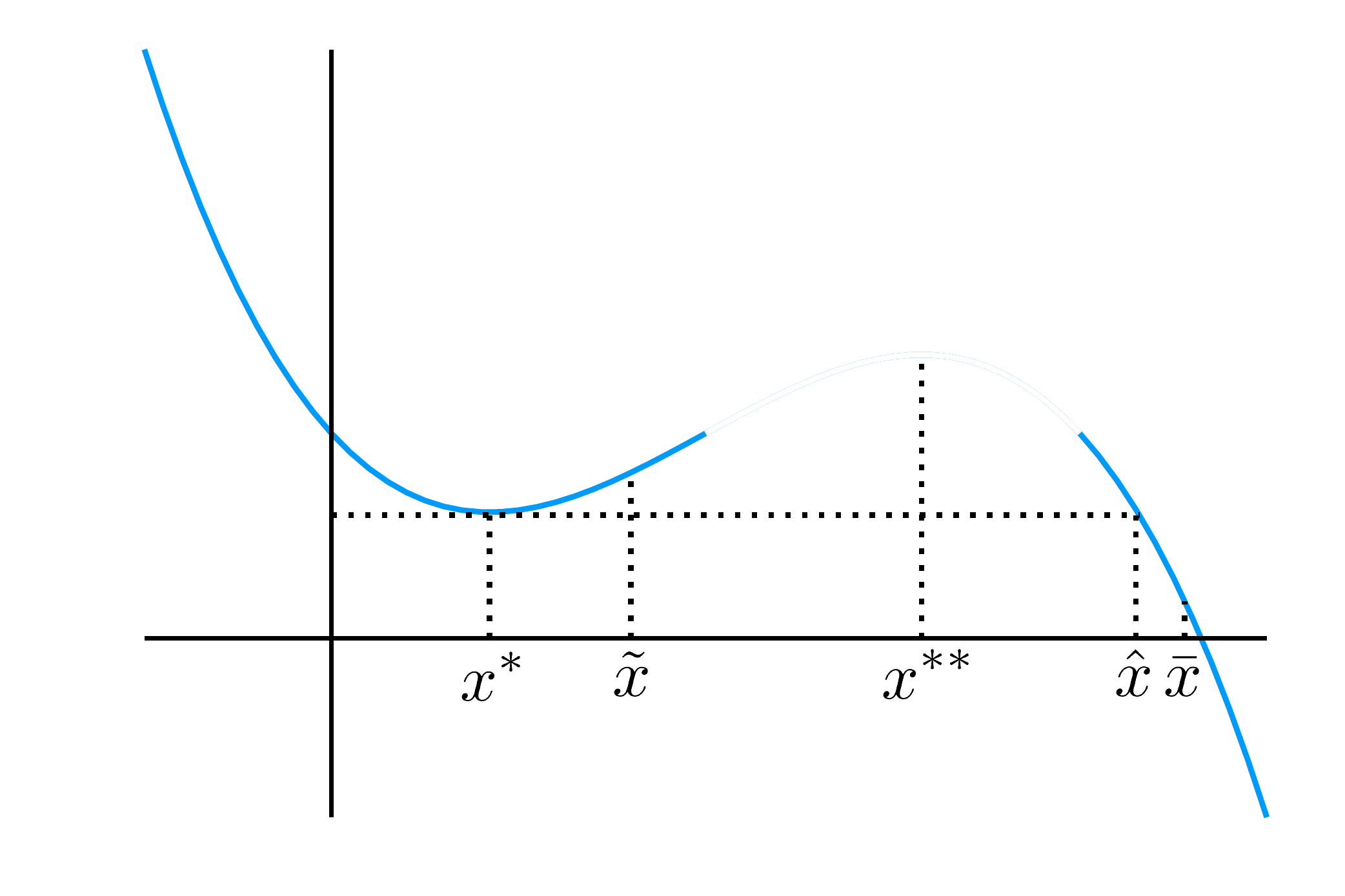}
\caption{Illustration of the proof of Theorem \ref{th:single_global}.}
\label{fig_single_global}
\end{figure}

\begin{theorem}
\label{th:single_global}
Let $f:(a,b)\subset\R\to\R$ be a differentiable function with a single critical point $x^*\in(a,b)$. 
If $x^*$ is a local minimizer, then it is a global minimizer. 
\end{theorem}
\beginproof
Assume by contradiction that there exists $\bar x\in(a,b)$, say, $\bar x>x^*$, such that 
$f(\bar x)<f(x^*)$. Since $x^*$ is a local minimizer, there exists $\delta>0$ such that 
$f(x^*)\leq f(x)$ for all $x\in(x^*-\delta,x^*+\delta)$. In fact, we have $f(x^*)<f(x)$ for all 
$x\in(x^*-\delta,x^*+\delta)\setminus\{x^*\}$, because otherwise there would be another critical point 
of $f$, in view of Rolle's theorem. Consider then $\tilde x\in(x^*,\bar x)$ with $f(x^*)<f(\tilde x)$. 
So, the intermediate value theorem guarantees the existence of $\hat x\in(\tilde x,\bar x)$ with 
$f(\hat x)=f(x^*)$. Therefore, by the Rolle's theorem, we conclude that exists a critical point 
$x^{**}\in(x^*,\hat x)$, contradicting the hypothesis. Figure \ref{fig_single_global} illustrates this 
proof.
\endproof

It is well known that the converse of Theorem \ref{th:kkt} is not necessarily true. That is, 
the optimality conditions (\ref{kkt_grad})--(\ref{kkt_feas}) are not sufficient to ensure 
that the point is a local minimizer. Indeed, these conditions are also satisfied at a 
maximizer.
When dealing with unconstrained minimization in two variables, we have the famous sufficient 
condition, present in almost all textbooks on the subject, to ensure that a critical point 
$x^*$ is a local minimizer, namely, the test of second derivatives
$$
\dfrac{\partial^2f}{\partial x_1^2}(x^*)>0\quad\mbox{and}\quad
\dfrac{\partial^2f}{\partial x_1^2}(x^*)\dfrac{\partial^2f}{\partial x_2^2}(x^*)-
\left(\dfrac{\partial^2f}{\partial x_1\partial x_2}(x^*)\right)^2>0.
$$
However, it is not so common to discuss a test for constrained optimization. In the next 
remark we address another issue related to this subject. 

\begin{remark}
\label{rmk:ctex2}
{\em In the context of problem (\ref{eq_prob}), it is also typical the following
question: 
{\bf ``Among the critical points, is the one with the smallest 
image (under the function) a local minimizer?".}
Again, the answer is no and the example below shows why.}
\end{remark}

\begin{figure}[htbp]
\centering
\includegraphics[scale=0.45]{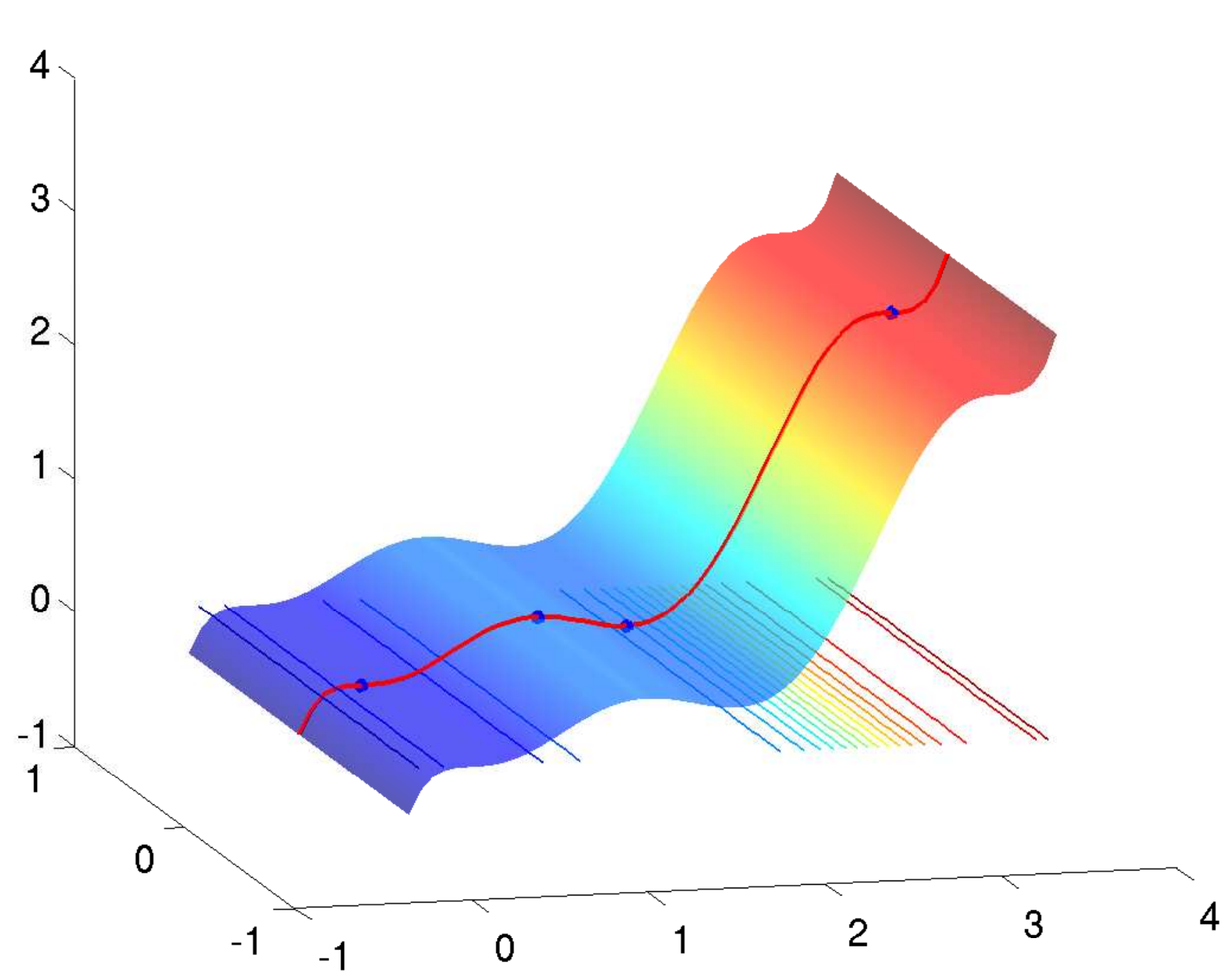}
\includegraphics[scale=0.38]{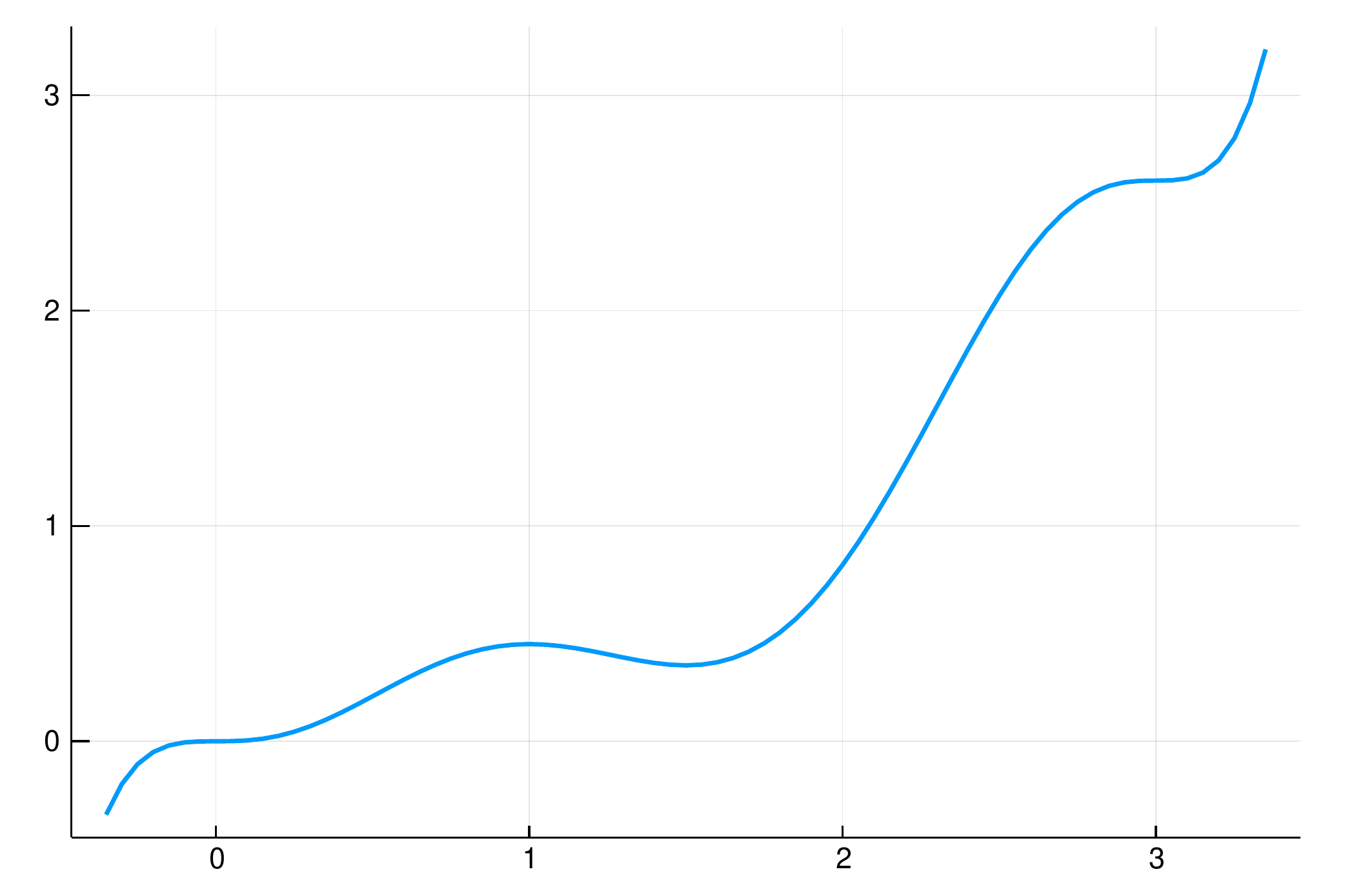}
\caption{Illustration of Example \ref{ex:small}. The left picture shows the graph 
of $f$ and the curve on it corresponding to the feasible set. On the right, we have 
the graph of the auxiliary function $\varphi$, defined in (\ref{auxphi}). 
Note that there are four critical points, being two of them local extrema. Besides, 
the smallest critical value does not correspond to a local minimizer.}
\label{fig_ctex2}
\end{figure}

\begin{example}
\label{ex:small}
Consider the problem (\ref{eq_prob}) with the functions $f,g:\R^2\to\R$ defined by 
$$
f(x)=\dfrac{1}{7}x_1^7-\dfrac{17}{12}x_1^6+
\dfrac{51}{10}x_1^5-\dfrac{63}{8}x_1^4+\dfrac{9}{2}x_1^3
$$ 
and $g(x)=x_2$. Find the critical points, its images and say which one is a minimizer. 
\end{example}
\beginsol
The condition (\ref{kkt_grad}) in this case is 
$$
\dfrac{1}{2}\left(\begin{array}{c} x_1^2(x_1-3)^2(x_1-1)(2x_1-3) \\ 
0 \end{array}\right)+\lambda\left(\begin{array}{c} 0 \\ 1 \end{array}\right)=
\left(\begin{array}{c} 0 \\ 0 \end{array}\right),
$$
yielding $\lambda^*=0$, $\bar x_1=0$, $\hat x_1=1$, $x_1^*=\frac{3}{2}$ and 
$\tilde x_1=3$. So, we have four critical points 
$$
\bar x=\left(\begin{array}{r} 0 \\ 0 \end{array}\right),\;\;
\hat x=\left(\begin{array}{r} 1 \\ 0 \end{array}\right),\;\;
x^*=\dfrac{3}{2}\left(\begin{array}{r} 1 \\ 0 \end{array}\right)\quad\mbox{and}\quad 
\tilde x=\left(\begin{array}{r} 3 \\ 0 \end{array}\right).
$$
By restricting the objective function to the feasible set, that is, to the points of the 
form $x=\left(\begin{array}{r} t \\ 0 \end{array}\right)$, with $t\in\R$, we obtain 
$f(x)=\varphi(t)$, where 
\begin{equation}
\label{auxphi}
\varphi(t)=\dfrac{1}{7}t^7-\dfrac{17}{12}t^6+\dfrac{51}{10}t^5-\dfrac{63}{8}t^4+
\dfrac{9}{2}t^3.
\end{equation}
Since $\varphi'(t)=\frac{1}{2}t^2(t-3)^2(t-1)(2t-3)$, we conclude that $\bar x$ is 
neither maximizer nor minimizer, but a saddle point. The same is true 
for $\tilde x$. On the other hand, $\hat x$ is a local maximizer and $x^*$ is a 
local minimizer. Finally, the critical values are  
$$
f(\bar x)=0,\;\;f(\hat x)\approx 0.4511,\;\;f(x^*)\approx 0.3525
\quad\mbox{and}\quad f(\tilde x)\approx 2.6035.
$$
It should be noted that the smallest critical value does not correspond to a local 
minimizer and that the greatest critical value does not correspond to a local maximizer.
Figure \ref{fig_ctex2} illustrates this example.
\endsol

The next section is devoted to discuss sufficient conditions to ensure optimality for constrained 
optimization problems. 

\section{Sufficient optimality conditions}
\label{sec_suff}
In this section we present a criterion, based on the second derivatives, for attesting that 
a critical point is a local minimizer. For completeness we present first a general result, well 
known in the optimization community. Then, we particularize the test to the specific cases 
studied in Calculus.

We stress that the criteria we will present are only local conditions and do not say anything 
about global minimization without additional assumptions or specific situations. 

To simplify the presentation consider the Lagrangian function associated with the 
problem (\ref{eq_prob}),
$$
(x,\lambda)\in\mathbb{R}^n\times\mathbb{R}^m\mapsto \ell(x,\lambda)=f(x)+\lambda^T g(x).
$$
The Lagrangian Hessian, that is, the matrix of the second derivatives of $\ell$ with 
respect to $x$, is denoted by
$$
\nabla_{xx}^2\ell(x,\lambda)=\nabla^2 f(x)+\sum_{i=1}^{m}\lambda_i\nabla^2 g_i(x).
$$
The result below can be found in many optimization books. See, for 
example, \cite{luenberger,nocedal}.

\begin{theorem}
\label{th:general}
Let $x^*\in\R^n$ be a critical point for the problem (\ref{eq_prob}) and $\lambda^*\in\R^m$ a 
corresponding multiplier vector, according to Definition \ref{def:stat}. Suppose that
$$
d^T\nabla_{xx}^2\ell( x^*,\lambda^*)d>0
$$
for all nonzero vectors $d\in\R^n$ satisfying $\nabla g_i(x^*)^Td=0$, $i=1,\ldots,m$. Then there 
exist $\delta>0$ and a neighborhood $V$ of $ x^*$ such that
$$
f(x)-f( x^*)\geq\delta\|x- x^*\|^2
$$
for all $x\in V$ with $g(x)=0$. In particular, $ x^*$ is a strict local minimizer of (\ref{eq_prob}).
\end{theorem}

Despite the existence of this condition for general dimensions, we consider here the 
particular $2$ and $3$-dimensional cases with one or two constraints, which are the most 
common cases in the Calculus courses. In these situations, the Hessian matrices of a 
function $\varphi$ are 
$$
\nabla^2\varphi=\left(\begin{array}{cc} 
\dfrac{\partial^2\varphi}{\partial x_1^2} & \dfrac{\partial^2\varphi}{\partial x_1\partial x_2} 
\vspace{5pt} \\ 
\dfrac{\partial^2\varphi}{\partial x_2\partial x_1} & \dfrac{\partial^2\varphi}{\partial x_2^2}
\end{array}\right)
\quad\mbox{or}\quad
\nabla^2\varphi=\left(\begin{array}{ccc} 
\dfrac{\partial^2\varphi}{\partial x_1^2} & \dfrac{\partial^2\varphi}{\partial x_1\partial x_2} 
 & \dfrac{\partial^2\varphi}{\partial x_1\partial x_3} \vspace{5pt} \\ 
\dfrac{\partial^2\varphi}{\partial x_2\partial x_1} & \dfrac{\partial^2\varphi}{\partial x_2^2} 
& \dfrac{\partial^2\varphi}{\partial x_2x_3} \vspace{5pt} \\ 
\dfrac{\partial^2\varphi}{\partial x_3\partial x_1} & \dfrac{\partial^2\varphi}{\partial x_3x_2} 
& \dfrac{\partial^2\varphi}{\partial x_3^2}
\end{array}\right)
$$
if $n=2$ or $n=3$, respectively.

\subsection{The two variables and one constraint case}
\label{sec_suff_R2}
Consider the problem (\ref{eq_prob}) with $n=2$ and $m=1$. That is, the problem
of minimizing a function of two variables subject to a single equality 
constraint.\\
\noindent The next theorem follows immediately from the previous one.

\begin{theorem}
\label{th:kkt2a}
Let $x^*\in\R^2$ be a critical point for the problem (\ref{eq_prob}) and let $\lambda^*\in\R$ 
be the corresponding Lagrange multiplier. Define $H=\nabla^2f(x^*)+\lambda^*\nabla^2g(x^*)$, 
assume that $\nabla g(x^*)\neq 0$ and take a nonzero vector $v\perp\nabla g(x^*)$. If 
\begin{equation}
\label{eq_kkt2a}
v^THv>0,
\end{equation}
then $x^*$ is a local minimizer for the problem. 
\end{theorem}

Now, let us see a straighforward application of the previous theorem.

\begin{figure}[htbp]
\centering
\includegraphics[scale=0.415]{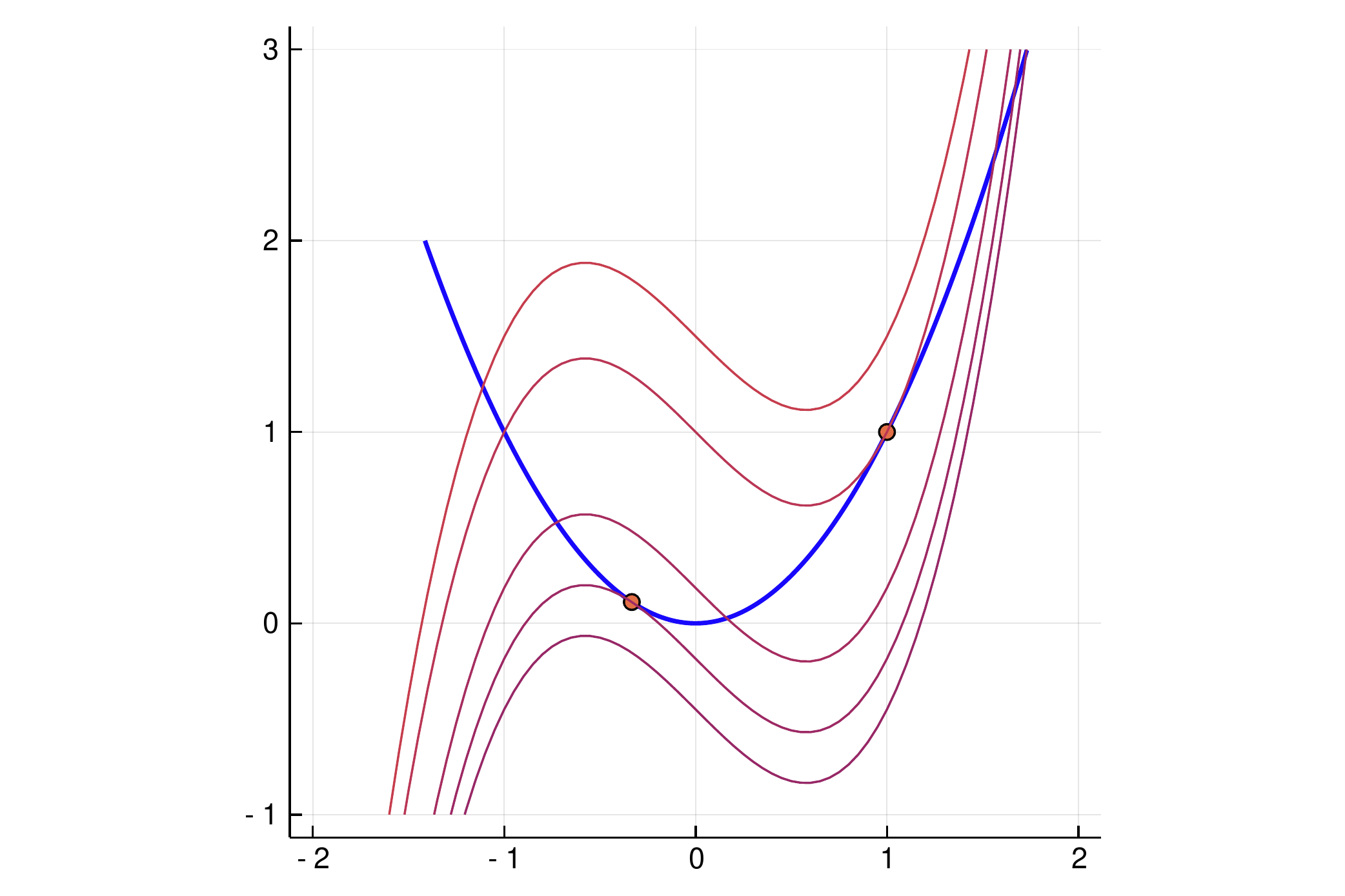}
\includegraphics[scale=0.42]{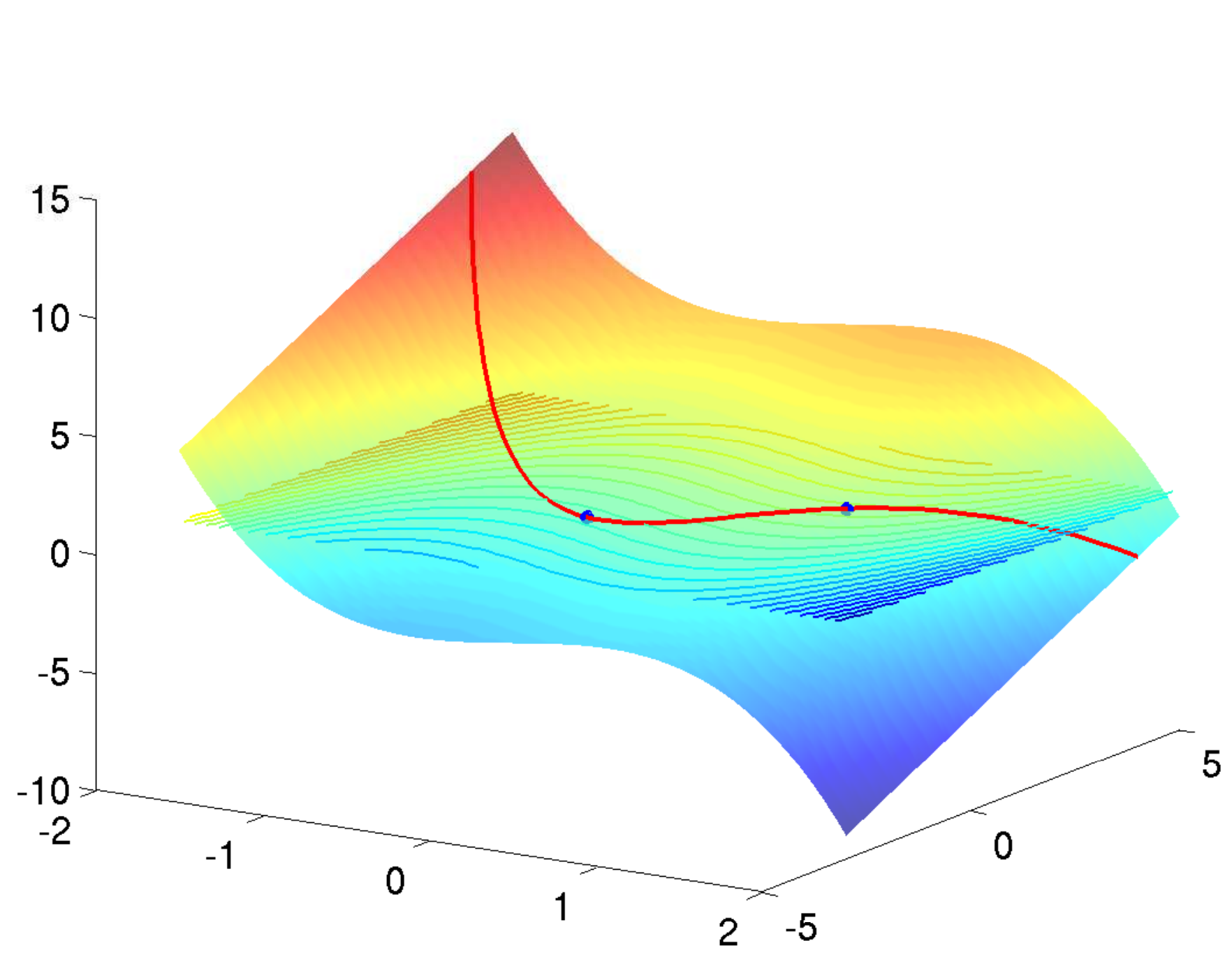}
\caption{On the left, the level curves of $f$, given by $x_2=x_1^3-x_1+constant$,  
the feasible set, given by $x_2=x_1^2$, and the two critical points $x^*$ and $\bar x$. 
On the right, the graph of $f$ and the points on it corresponding to the feasible set.}
\label{fig_exsuff1}
\end{figure}

\begin{example}
\label{ex:suff1}
Discuss the problem (\ref{eq_prob}) with $f,g:\R^2\to\R$ defined by $f(x)=x_2-x_1^3+x_1$ and 
$g(x)=x_2-x_1^2$.
\end{example}
\beginsol
The condition (\ref{kkt_grad}) in this case is 
$$
\left(\begin{array}{c} -3x_1^2+1 \\ 1 \end{array}\right)+\lambda
\left(\begin{array}{c} -2x_1 \\ 1 \end{array}\right)=
\left(\begin{array}{c} 0 \\ 0 \end{array}\right),
$$
giving $\lambda^*=-1$ and $x_1=-\dfrac{1}{3}$ or $x_1=1$. So, we have two critical points 
$$
x^*=\dfrac{1}{9}\left(\begin{array}{r} -3 \\ 1 \end{array}\right)\quad\mbox{and}
\quad \bar x=\left(\begin{array}{r} 1 \\ 1 \end{array}\right).
$$
Moreover, 
$$
H=\left(\begin{array}{cc} 2-6x_1 & 0 \\ 0 & 0 \end{array}\right)\quad\mbox{and}
\quad v=\left(\begin{array}{c} 1 \\ 2x_1 \end{array}\right)\perp\nabla g(x).
$$
Thus, $x^*$ is a local minimizer since $v^THv=2-6x_1^*>0$. Note that this point is not a 
global minimizer because 
$f\left(\begin{array}{r} 2 \\ 4 \end{array}\right)=-2<-\dfrac{5}{27}=f(x^*)$. 
Figure \ref{fig_exsuff1} illustrates this example.
\endsol

\begin{remark}
\label{rmk:suff1}
{\em In the context of Theorem \ref{th:kkt2a} we also have the maximization condition. If 
$v^THv<0$, then $x^*$ is a local maximizer of $f(x)$ subject to $g(x)=0$.}
\end{remark}

\subsection{The three variables and one constraint case}
\label{sec_suff_R3}
Consider the problem of minimizing a function of three variables subject to a 
single equality constraint.\\ 
\noindent Here comes another simple application of the general theorem:

\begin{theorem}
\label{th:kkt3a}
Let $x^*\in\R^3$ be a critical point for the problem (\ref{eq_prob}) with $n=3$ and $m=1$. 
Consider $\lambda^*\in\R$ the corresponding Lagrange multiplier and define 
$H=\nabla^2f(x^*)+\lambda^*\nabla^2g(x^*)$. Suppose that $\nabla g(x^*)\neq 0$ and take 
vectors $v_1,v_2\in\R^3$ such that ${\rm span}\{v_1,v_2\}=\nabla g(x^*)^\perp$. Consider 
the matrices $V=(v_1\ v_2)\in\R^{3\times 2}$ and $A=V^THV\in\R^{2\times 2}$. If 
$$
a_{11}>0\quad\mbox{and}\quad\det(A)>0,
$$
then $x^*$ is a local minimizer for the problem. 
\end{theorem}

Let us see now a straighforward application of the previous theorem.

\begin{example}
\label{ex:suff3}
Let us revisit Example \ref{ex:min_area}. Suppose we have not proved the existence of 
a global minimizer. Then, we are able to establish (only) the local condition. 
\end{example}
\beginsol
We have the critical point 
$$
x^*=\dfrac{1}{2}\left(\begin{array}{c} 2\sqrt[3]{2} \\ 2\sqrt[3]{2} \\ \sqrt[3]{2}
\end{array}\right)
$$
with multiplier $\lambda^*=-2\sqrt[3]{4}$. Thus, 
$$
H=-\left(\begin{array}{ccc} 0 & 1 & 2 \\ 1 & 0 & 2 \\ 2 & 2 & 0 \end{array}\right)
\quad\mbox{and}\quad \nabla g(x^*)=\dfrac{\sqrt[3]{4}}{2}
\left(\begin{array}{c} 1 \\ 1 \\ 2 \end{array}\right).
$$
We can consider, for example, 
$
\nabla g(x^*)^\perp={\rm span}\left\{
\left(\begin{array}{r} 1 \\ -1 \\ 0 \end{array}\right),
\left(\begin{array}{r} 2 \\ 0 \\ -1 \end{array}\right)
\right\},
$
yielding 
$$
A=-\left(\begin{array}{crr} 1 & -1 & 0 \\ 2 & 0 & -1 \end{array}\right)
\left(\begin{array}{ccc} 0 & 1 & 2 \\ 1 & 0 & 2 \\ 2 & 2 & 0 \end{array}\right)
\left(\begin{array}{rr} 1 & 2 \\ -1 & 0 \\ 0 & -1 \end{array}\right)
=\left(\begin{array}{cc} 2 & 2 \\ 2 & 8 \end{array}\right)
$$
Hence, $x^*$ is a local minimizer since $a_{11}=2>0$ and $\det(A)=12>0$. 
\endsol

\subsection{The three variables and two constraints case}
\label{sec_suff_R3g2}
Consider the problem of minimizing a function of three variables subject to two equality 
constraints. 

Finally, let us present our last application of the general theorem.

\begin{theorem}
\label{th:kkt3b}
Let $x^*\in\R^3$ be a critical point for the problem (\ref{eq_prob}) with $n=3$ and $m=2$. 
Consider $\lambda^*\in\R^2$ the vector of Lagrange multipliers and define 
$H=\nabla^2f(x^*)+\lambda_1^*\nabla^2g_1(x^*)+\lambda_2^*\nabla^2g_2(x^*)$. Suppose that 
$\nabla g_1(x^*)$ and $\nabla g_2(x^*)$ are linearly independent and take a nonzero vector 
$v\in\R^3$ such that $v\perp\nabla g_1(x^*)$ and $v\perp\nabla g_2(x^*)$. If 
\begin{equation}
\label{eq_kkt3b}
v^THv>0,
\end{equation}
then $x^*$ is a local minimizer for the problem. 
\end{theorem}

\begin{figure}[htbp]
\centering
\includegraphics[scale=0.42]{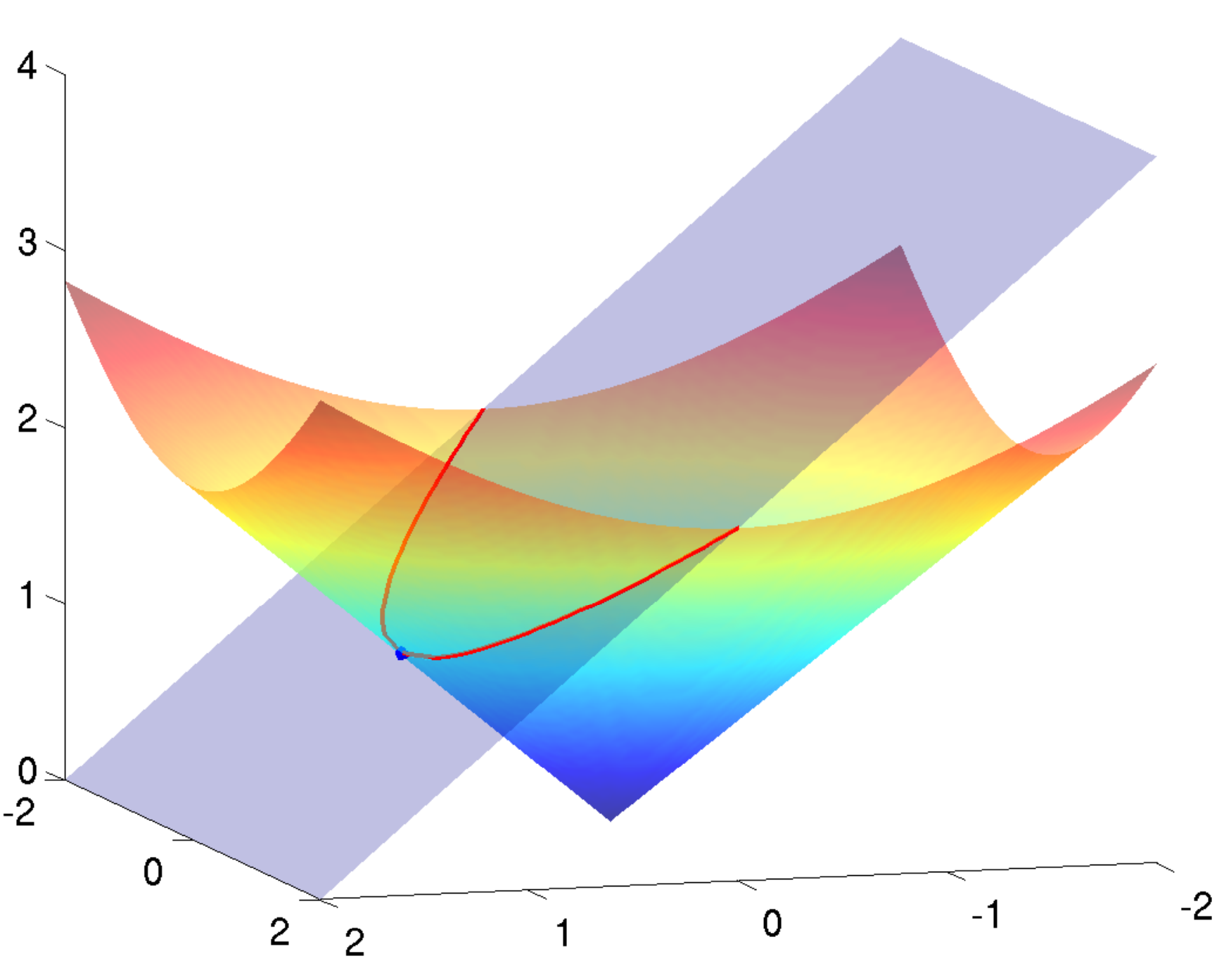}
\caption{Illustration of Example \ref{ex:suff3b}: the cone given by $x_1^2+x_2^2-x_3^2=0$, 
the plane defined by $x_1+x_3-2=0$ and the feasible set, the intersection represented 
by the red curve. In this problem we want to find the lowest point of that curve.}
\label{fig_exsuff3}
\end{figure}

\noindent Here comes our last example:

\begin{example}
\label{ex:suff3b}
Consider $f,g_1,g_2:\R^3\to\R$ defined by $f(x)=x_3$, $g_1(x)=x_1^2+x_2^2-x_3^2$ and 
$g_2(x)=x_1+x_3-2$. Solve the problem (\ref{eq_prob}) for these functions.
\end{example}
\beginsol
The condition (\ref{kkt_grad}) in this case is 
$$
\left(\begin{array}{c} 0 \\ 0 \\ 1 \end{array}\right)+\lambda
\left(\begin{array}{r} 2x_1 \\ 2x_2 \\ -2x_3 \end{array}\right)+\mu
\left(\begin{array}{c} 1 \\ 0 \\ 1 \end{array}\right)=
\left(\begin{array}{c} 0 \\ 0 \\ 0 \end{array}\right),
$$
which immediately implies that $\lambda^*\neq 0$ and hence, $x_2^*=0$. Using the constraints, 
we conclude that $x_1^*=x_3^*=1$. This in turn implies that $\lambda^*=\dfrac{1}{4}$ and 
$\mu^*=-\dfrac{1}{2}$. So, 
$$
H=\dfrac{1}{4}\left(\begin{array}{ccr} 2 & 0 & 0 \\ 0 & 2 & 0 \\ 0 & 0 & -2 \end{array}\right)\ ,\ 
g_1(x^*)=\left(\begin{array}{r} 2 \\ 0 \\ -2 \end{array}\right)
\quad\mbox{and}\quad \nabla g_2(x^*)=\left(\begin{array}{c} 1 \\ 0 \\ 1 \end{array}\right).
$$
We can consider, for example, 
$$
v=\left(\begin{array}{c} 0 \\ 1 \\ 0 \end{array}\right)\in
\left\{\nabla g_1(x^*),\nabla g_2(x^*)\right\}^\perp
$$
and see that $v^THv>0$, proving then that $x^*$ is a local minimizer for the problem. 
In fact, we can prove that this point is a global minimizer. To see this, note that 
$$
x_3^2=x_1^2+x_2^2=(2-x_3)^2+x_2^2=4-4x_3+x_3^2+x_2^2,
$$
which gives $4(x_3-1)=x_2^2\geq 0$. So, any feasible point $x$ satisfies 
$f(x)=x_3\geq 1=f(x^*)$.
Figure \ref{fig_exsuff3} illustrates the feasible set of this example.
\endsol

\begin{remark}
\label{rmk:v}
{\em It is easy to see that the condition (\ref{eq_kkt3b}) does not depend on the particular 
choice of $v$: if $v^THv>0$ for a vector $v\in\left\{\nabla g_1(x^*),\nabla g_2(x^*)\right\}^\perp$, 
then $\bar v^TH\bar v>0$ for any other nonzero vector 
$\bar v\in\left\{\nabla g_1(x^*),\nabla g_2(x^*)\right\}^\perp$. Indeed, in this case, 
$\bar v=\alpha v$, for some $\alpha\in\R\setminus\{0\}$. The same reasoning is true for 
condition (\ref{eq_kkt2a}). It can be also proved that the conditions in Theorem \ref{th:kkt3a} 
do not depend on the choice of vectors $v_1,v_2$ such that 
${\rm span}\{v_1,v_2\}=\nabla g(x^*)^\perp$.}
\end{remark}

\section{Conclusion}
\label{sec_concl} 
In this paper, we have pointed out that, in some examples (of some 
undergraduate Calculus textbooks) related to the acquirement of global 
minimizers via the Lagrange Multiplier Method (LMM), a little bit of imprecision
has been typical, particularly when dealing with worked problems. One way to 
mitigate that would be the use of a criterion to guarantee when a critical point 
(obtained by the LMM) is a local minimizer. So, we have proposed such a 
criterion, which, by the way, has been kept absent from Calculus textbooks.
On the other hand, for those Professors who jump into the `global' aspects of
the LMM, in spite of being a strictly local result, based on what we discussed
here, we also propose the following way to state the LMM:\\
\noindent{\it For continuously differentiable functions, $f$ and $g$, in order 
to determine the minimum value of $f$ subject to the constraint $g=k$ with $k$ 
constant, assuming that this global minimum value is attained on the interior of
the domain shared by $f$ and $g$, but not on the boundary of it, and that 
$\nabla g\neq\vec{0}$ holds for that domain, do the following:
\begin{enumerate}
\item Determine each point and, if necessary, also $\lambda$, satisfying the 
following system:
\begin{enumerate}
\item $\nabla f=\lambda\nabla g$;
\item g=k.
\end{enumerate}
\item Evaluate $f$ for those points obtained in the previous item: the smallest
value of $f$ is its global minimum.
\end{enumerate}}


\clearpage

\end{document}